\documentclass[12pt]{article}
\usepackage{amsfonts,amssymb,amsmath,graphics,theorem}
\def\lape{\kappa}

\newtheorem{Lem}{Lemma} 
\newtheorem{Thm}{Theorem}

\newtheorem{Cor}{Corollary}

\newtheorem{Pro}{Proposition}

\def\be{\begin{enumerate}}  \def\bi{\begin{itemize}}
\def\ee{\end{enumerate}}   \def\ei{\end{itemize}} \def\im{\item}
\def\ssec{\subsection}   \def\sec{\section} 
 
\def\no{\nonumber}

\def\nn{\nonumber}
\def\beq{\begin{eqnarray}} \def\eeq{\end{eqnarray}}

\def\al*#1{\begin{align*}#1\end{align*}}

\def\ga*#1{\begin{gather*}#1\end{gather*}}

\def\alat*#1#2{\begin{alignat*}{#1}#2\end{alignat*}}
\def\bea{\begin{eqnarray*}}
\def\eea{\end{eqnarray*}}
\def\ml*#1{\begin{multline*}#1\end{multline*}}

 \def\mbf{\mathbf} 
  \def\ovl{\overline}
\newcommand{\Bf}[1]{{\mbox{\scriptsize\boldmath$#1$}}}
\newcommand{\bff}[1]{{\mbox{\boldmath$#1$}}}
\def\bcm{\begin{comment}} \def\ecm{\end{comment}}

 \def\I{\int}

\def\mc{\mathcal}

 \def\le{\left} \def\ri{\right}  \def\i{\infty}

\def\te#1{\mathrm{e}^{#1}}  \def\td{\text{\rm d}}

\def\T{\tilde}  
 \def\WT{\widetilde}

\def\a{\alpha} \def\b{\beta}
\def\g{\gamma}  \def\d{\delta}   \def\th{\theta}

\def\e{\epsilon} \def\k{\kappa} \def\l{\lambda} \def\m{\mu} 
    \def\r{\rho} \def\s{\sigma}
\def\t{\tau}     \def\ps{\psi}
 \def\w{\omega} \def\q{\qquad} 
 \def\G{\Gamma}

\newcommand{\halmos}{{\mbox{\, \vspace{3mm}}} \hfill\mbox{$\square$}}

\def\la{\label}  
 \def\le{\left} \def\ri{\right}

 \def\srui{\ovl{\psi}}

    \def\lape{\kappa}

\begin{document}

\def\theequation{\arabic{equation}} % If you don't want the equations counter

                                     % to be reset with each section,

                                     % outcomment this definition

\title{A two-dimensional ruin problem on the positive quadrant.}

\author{Florin Avram
\footnote{ Dept. de Math., Universit\'e de Pau, E-mail:
Florin.Avram@univ-Pau.fr} \q Zbigniew
Palmowski\footnote{University of Wroclaw, pl. Grunwaldzki 2/4,
50-384 Wroclaw, Poland and Utrecht University, P.O. Box 80.010,
3500 TA, Utrecht, The Netherlands, E-mail:
zpalma@math.uni.wroc.pl} \q and Martijn
Pistorius\footnote{Department of Mathematics, King's College
London, Strand, London WC2R 2LS, UK, Email:
Martijn.Pistorius@kcl.ac.uk } }
\date{}

\maketitle

%\begin{center}\fbox{{\large Working %Notes}}\\\end{center}

 \begin{abstract}
In this paper we study the joint ruin problem for two insurance
companies that divide between them both claims  and premia in some
specified proportions (modeling two branches of the same insurance
company or an insurance and re-insurance company). Modeling the risk
processes of the insurance companies by Cram\'{e}r-Lundberg
processes we obtain the Laplace transform in space of the
probability that either of the insurance companies is ruined in
finite time. Subsequently, for exponentially distributed claims, we
derive an explicit analytical expression for this joint ruin
probability by explicitly inverting this Laplace transform.
We also provide a characterization of the Laplace transform
of the joint ruin time.
\end{abstract}

\newpage

\sec{A  two dimensional ruin problem \label{s:intr}} In this paper
we consider a particular two dimensional risk model in which two
companies split the amount they pay out of each claim in positive
proportions $\d_1$ and $\d_2$ with $\d_1+\d_2=1$,
and the premiums according to rates $c_1$ and $c_2$. Thus, the risk
process $U_i$ of the $i$'th company satisfies
\[U_i(t):=-\delta_i S(t)+c_i t+u_i,\qquad i=1,2\; ,\]
where $u_i$ are the initial reserves. We will work with a spectrally
positive L\'{e}vy process $S(t)$, that is L\'{e}vy process with only
upward jumps that represents the cumulative amount of claims up to
time $t$. In particular we focus on the classical
Cram\'{e}r-Lundberg model:
\begin{equation}\label{St} S(t)=\sum_{k=1}^{N(t)}
\sigma_k,\end{equation} where $N(t)$ is a Poisson process with
intensity $\lambda$ and the claims $\sigma_k$ are i.i.d. random
variables independent of $N(t)$, with distribution function $F(x)$
and mean  $E[\sigma_k]=\m^{-1}$. We shall assume that the second
company, to be called the reinsurer, gets smaller profits per amount
paid, i.e.: \beq p_1=\frac{ c_1}{\delta_1} >
 \frac{ c_2}{\delta_2}=p_2. \label{inter}\eeq
As usual in risk theory, we  assume that
$p_i>\rho:=\frac{\lambda}{\m}$, which implies that { in the
absence of ruin}, $U_i(t)\to\infty$ as $t\to\infty$ ($i=1,2$).
Ruin happens  at the time $\tau=\tau(u_1,u_2)$ when at least one
insurance company is  ruined:
\begin{equation}\label{eq:ruint}
\t(u_1,u_2) := \inf\{t\ge 0: U_1(t) <  0 \quad\mbox{or}\quad U_2(t)
<  0\},
\end{equation}
i.e. at the first exit time of $(U_1(t), U_2(t))$ from the positive
quadrant. In this paper we will analyze the perpetual or ultimate
ruin probability:
\begin{equation}\label{eq:ruin}
\psi(u_1,u_2) = P\left[\t(u_1,u_2) <\infty\right].
\end{equation}
 Although ruin theory
under multi-dimensional models rarely admits analytical solutions,
we are able to obtain in our problem a closed form solution for
(\ref{eq:ruin}) if $\sigma_i$ are exponentially distributed with
intensity $\mu$.
\bigskip

\begin{figure}[t]
\begin{center}
\leavevmode \resizebox{0.62 \textwidth}{!}{\setlength{\unitlength}{2000sp
%3947sp
}%
\begingroup\makeatletter\ifx\SetFigFont\undefined%
\gdef\SetFigFont#1#2#3#4#5{%
  \reset@font\fontsize{#1}{#2pt}%
  \fontfamily{#3}\fontseries{#4}\fontshape{#5}%
  \selectfont}%
\fi\endgroup%
\begin{picture}(10062,7362)(1114,-7711)
\thinlines {\put(1201,-361){\line( 0,-1){7200}}
\put(1201,-7561){\line( 1, 0){9600}}
}%
{\put(1201,-361){\line( -1, -2){ 175}}
}%
{\put(1201,-361){\line( 1,-2){ 175}}
}%
{\put(10450,-7400){\line( 2,-1){350}}
}%
{\put(10450,-7736){\line( 2, 1){350}}
}%
{\multiput(1201,-7561)(4.99383,7.49075){812}{\makebox(1.6667,11.6667){\SetFigFont{12}{14.4}{\rmdefault}{\mddefault}{\updefault}.}}
}%
{\put(2476,-5311){\line( 0,-1){375}} \put(2476,-5686){\line(
4,3){375}}
}%
{\put(1201,-7561){\line( 3, 1){8797.500}}
}%
{\put(4651,-6200){\line( 5,-1){450}}
\put(5126,-6275){\line(-3,-4){325}}
}%
\put(5626,-6511){\makebox(0,0)[lb]{\smash{\SetFigFont{12}{14.4}{\rmdefault}{\mddefault}{\updefault}{drift direction}%
}}}
\put(6701,-4086){\makebox(0,0)[lb]{\smash{\SetFigFont{12}{14.4}{\rmdefault}{\mddefault}{\updefault}{$(\delta_1  ,\delta_2)$}%
}}}
\put(11026,-7636){\makebox(0,0)[lb]{\smash{\SetFigFont{12}{14.4}{\rmdefault}{\mddefault}{\updefault}{$u_1$}%
}}}
\put(1501,-736){\makebox(0,0)[lb]{\smash{\SetFigFont{12}{14.4}{\rmdefault}{\mddefault}{\updefault}{$u_2$}%
}}}
\put(3751,-4111){\makebox(0,0)[lb]{\smash{\SetFigFont{12}{14.4}{\rmdefault}{\mddefault}{\updefault}{jump directions}%
}}}
\put(8476,-6511){\makebox(0,0)[lb]{\smash{\SetFigFont{12}{14.4}{\rmdefault}{\mddefault}{\updefault}{$(c_1,c_2)$}%
}}}
\end{picture}}

\caption{Geometrical considerations}\label{fig1}
\end{center}
\end{figure}
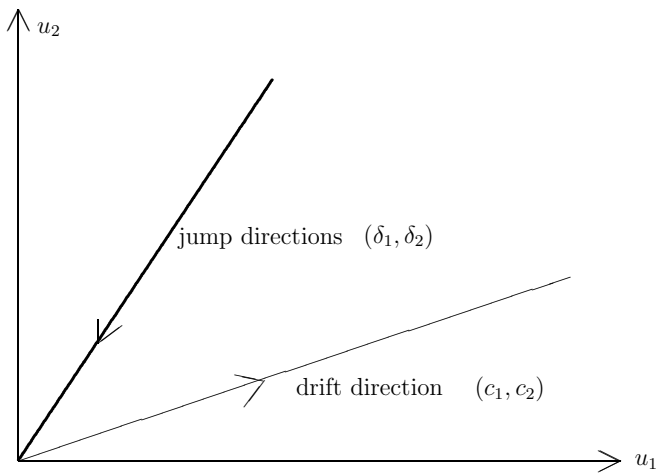

{\bf Geometrical considerations.} The solution of the
two-dimensional ruin problem (\ref{eq:ruin}) strongly depends on the
relative sizes of the proportions $\bff\d = (\d_1, \d_2)$ and
premium rates $\bff c = (c_1, c_2)$ -- see Figure \ref{fig1}.
If, as assumed throughout, the angle of the vector $\bff\d$ with
the $u_1$ axis is larger than that of $\bff c$, i.e.  ${\delta_2
c_1 > \delta_1 c_2}$ we note that starting with initial capital
$(u_1,u_2)\in\mc C$ in the cone $\mathcal{C}=\{(u_1,u_2): u_2\leq
(\d_2/\d_1)u_1\}$ situated below the line $u_2=(\delta_2/\delta_1)
u_1$, the process $(U_1,U_2)$ ends up hitting at time $\t$ the
$u_1$ axis. Thus, in the domain $\mc C$ ruin occurs iff there is
ruin in the one-dimensional problem corresponding to the risk
process $U_2$.

{\bf One dimensional reduction.} A key observation  is that $\t$ in
(\ref{eq:ruint}) is also equal to
$$\tau(u_1,u_2)=\inf\{t\ge 0: S(t)>b(t)\},$$
where $b(t)=\min\{(u_1+c_1t)/\delta_1, (u_2+c_2t)/\delta_2\}$. The
two dimensional problem (\ref{eq:ruin}) may thus be also viewed as
a one dimensional crossing problem over a {piecewise linear
barrier}.

In the case that the initial reserves $u_1$ and $u_2$ are such
that $(u_1,u_2)\in\mc C$, that is, $u_2/\delta_2 \leq
u_1/\delta_1$, the barrier $b$ is linear,
$b(t)=(u_2+c_2t)/\delta_2$, the ruin happens always for the second
company. Thus, as we already observed, the problem (\ref{eq:ruin})
reduces in fact to the classical one-dimensional ultimate ruin
problem with premium $c_2$ and claims $\d_2 \s$, i.e.
$$\ps(u_1,u_2) =
\psi_2(u_2):= P(\t_2(u_2) < \i),$$ where $\t_2(u_2) = \inf\{t\ge0:
U_2(t) < 0\}$ and $\psi_2(u_2) $ is the ruin probability of $U_2$,
with $U_2(0)=u_2$. For the model (\ref{St}) the Pollaczeck-Khinchine
formula, well known from the theory of one-dimensional ruin (see
e.g. \cite{Rolski} or \cite{as2}), yields then an explicit series
solution for $\psi(u_1,u_2) = \psi_2(u_2)$ in the case of a
general claims distribution. For the phase-type claims $(\bff \b,
\bff B)$, i.e. with $P[\s>x] = \bff \b \te {\Bf B x} \bff 1$, the
ruin probability
 may be written in a simpler matrix exponential form:
\begin{eqnarray*}
\psi_2(u_2) &=&  \bff \eta  \te{\d_2^{-1}(\Bf B + \Bf b \Bf \eta) u_2}
\bff 1
\end{eqnarray*}
with $\bff \eta= \frac{\lambda}{p_2} \bff \b ( - \bff B)^{-1}$
(see for example (4) in \cite{AAU}), and
 in the case of exponential claim sizes with intensity $\mu$, it
 reduces to:
\begin{equation}\label{exponruin} \psi_2(u_2)= C_2e^{-(\gamma_2/\d_2) u_2}\;,
\end{equation}
where $\gamma_2=\mu -\lambda\d_2/c_2 = \mu - \lambda/p_2 $ and
$C_2= \frac{\lambda\d_2}{\mu c_2} = \frac{\lambda}{\mu p_2}$.

The rest of the paper is devoted to the analysis of
the opposite case, $u_2/\delta_2 >u_1/\delta_1$ and is organised as follows.
Section \ref{exa} is devoted to the Laplace transform of
the ruin time in the case $S$ is of the form (\ref{St})
with exponential jumps. Subsequently, in Sections \ref{s:pro} and \ref{s:LT}
we derive the Laplace double transform in $(u_1,u_2)$
of the ruin probability $\psi(u_1,u_2)$ if $S$ is a general
spectrally positive L\'{e}vy process.
Finally, in Section \ref{sec:LT3} this Laplace transform is
explicitly inverted, in the case of exponential claim sizes.

\sec{\label{exa} The differential system
 for the exponential claim sizes case }

In this section we provide a system of partial differential equations
for the Laplace transform
\begin{equation}\label{eq:psis}
\psi(u_1,u_2,s) := E[e^{-s\tau(u_1,u_2)}\mbf 1_{\{\tau(u_1,u_2)<\i\}}]
\end{equation}
of the ruin time $\tau(u_1,u_2)$ in (\ref{eq:ruint}) in the case
that $S$ is given by a compound Poisson process (\ref{St})
with intensity $\lambda$ and with claims sizes $\s_i$ that are
exponentially distributed with parameter $\mu$. The memoryless
property of interarrival times and claim sizes opens up the
possibility of embedding the (discontinuous) Cram\'{e}r-Lundberg
processes $(U_1,U_2)$ into a continuous Markov-modulated fluid
model. Informally, this is achieved by a transformation that
replaces the jumps of $(U_1,U_2)$ by a linear movement in the
direction $(-\d_1,-\d_2)$ of  duration equal to the size of the
jump, creating thereby a new {\em continuous} semi-Markovian model
$(\T{{U}}_1, \T U_2)$ called the {\em fluid embedding} of
$(U_1,U_2)$. As the process $(\T{{U}}_1, \T U_2)$ crosses
boundaries continuously and has exactly the same maxima and minima
as $(U_1,U_2)$, first passage problems may be easier to handle for
$(\T{{U}}_1, \T U_2)$ than for the original process $(U_1,U_2)$.
A formal construction of the fluid-embedding is given in the Appendix.
If we write $\T\t(u_1,u_2)=\inf\{t\ge0: \min\{\T U_1(t), \T U_2(t)\}< 0\}$
for the joint ruin time of $(\T U_1,\T U_2)$, it follows from the definition
of $(\T U_1,\T U_2)$ that $\t(u_1,u_2) = I(\T\t(u_1,u_2))$, where
$I(t)$ denotes the time up to time $t$ that $(\T U_1,\T U_2)$
was increasing. In particular,
\begin{equation}\label{eq:psis2}
\psi(u_1,u_2,s) = E[e^{-s I(\T\t(u_1,u_2))}\mbf 1_{\{\T\t(u_1,u_2)<\i\}}|Y_0=1].
\end{equation}
Setting $\phi(u_1,u_2,s)= E[e^{-s I(\T\t(u_1,u_2))}\mbf
1_{\{\T\t(u_1,u_2)<\i\}}|Y_0=-1]$ and writing
$\psi_{u_i}$ and $\phi_{u_i}$ for the partial derivative of $\psi$
and $\phi$ with respect to $u_i$ we have the following
characterization of $\phi$ and $\psi$.

\begin{Thm} For $\d_2 u_1 \leq \d_1 u_2$ it holds that
$( {\psi}(u_1,u_2,s),$
$\phi(u_1,u_2,s))^{\top}$ solves the Feynman-Kac system:
\begin{align*}  \left(\begin{array}{ll}c_1 & 0\\ 0 & -\delta_1 \end{array}\right)
\left(\begin{array}{l} {\psi}_{u_1}\\ \phi_{u_1 }
\end{array}\right)+ \left(\begin{array}{ll}c_2 & 0\\ 0 & -\delta_2
\end{array}\right) \left(\begin{array}{l} {\psi}_{u_2}\\ \phi_{u_2}
\end{array}\right)\\+ \left(\begin{array}{ll}-\lambda-s & \lambda\\
\mu & -\mu
\end{array}\right) \left(\begin{array}{l} {\psi}\\ \phi
\end{array}\right)=
\left(\begin{array}{l}0\\  0
\end{array}\right)
\label{maineq}
\end{align*}
with the boundary condition:
\begin{equation}\label{boundmain}
\left\{\begin{array}{ll}  {\psi}({u_1 },\frac{\d_2}{\d_1} u_1)&=
C_2 e^{-\g_2 \;
\frac{1}{\delta_1} u_1} \qquad\mbox{for all ${u_1 }\ge 0$,}\\
\phi(0 ,u_2)&=1 \qquad\mbox{for all ${u_2} \ge 0$.}
\end{array}
\right. \end{equation}
\end{Thm}
{\it Proof.} Conditioning on the first shared claim occurrence
epoch we obtain:
\begin{eqnarray}
\nn
\lefteqn{\psi(u_1,u_2,s)=e^{-\lambda h}e^{-s h} \psi(u_1+c_1
h,u_2+c_2 h, s)}\\
\nn
&&+ \int_0^h\lambda e^{-\lambda
t}\;dt\int_0^{\frac{u_1+c_1t}{\delta_1}\wedge
\frac{u_2+c_2t}{\delta_2}}\mu e^{-\mu z}e^{-s
t}\psi(u_1+c_1t-\delta_1z, u_2+c_2 t-\delta_2 z,s)dz\\
\label{eq:psiu1u2}
&&+\int_0^h\lambda e^{-\lambda
t}\;dt\int_{\frac{u_1+c_1t}{\delta_1}\wedge
\frac{u_2+c_2t}{\delta_2}}^\infty\mu e^{-\mu z}dz;\\
\lefteqn{\phi(u_1,u_2,s)=\int_0^\i\mu e^{-\mu z}\psi(u_1-\d_1 z,
u_2-\d_2 z, s)dz}. \label{eq:phiu1u2}
\end{eqnarray}
Both integrals on the LHS of (\ref{eq:psiu1u2})
go to 0 as $h\to 0$. This implies that
$\psi$ is a continuous function with respect to $u_1$ and $u_2$.
Note that $\frac{u_1+c_1t}{\delta_1}\wedge
\frac{u_2+c_2t}{\delta_2}= \frac{u_1+c_1t}{\delta_1}$ since we
live in the upper cone. This gives after simple manipulation:
\begin{eqnarray*}
\lefteqn{\frac{\psi(u_1+c_1h,u_2+c_2h,s)-\psi(u_1,u_2,s)}{h}+\frac{e^{-\lambda
h}e^{-s h}-1}{h} \psi(u_1+c_1 h,u_2+c_2 h, s)}\\&&+
\frac{1}{h}\int_0^h\lambda e^{-\lambda
t}\;dt\int_0^{(u_1+c_1t)/\delta_1}\mu e^{-\mu z}e^{-s
t}\psi(u_1+c_1t-\delta_1z, u_2+c_2 t-\delta_2 z,s)dz\\
&&+\frac{1}{h}\int_0^h\lambda e^{-\lambda t}\;dte^{-\mu
(u_1+c_1t)/\delta_1}=0.
\end{eqnarray*}
Note that the last 3 terms on the LHS of above equation have
limits as $h\to 0$, since $\psi$ is a continuous function. Thus
$\psi$ is a differentiable function of $(u_1,u_2)$.   Taking $h\to
0$ we derive
\begin{eqnarray*}
\lefteqn{c_1\psi_{u_1}(u_1,u_2,s)+c_2\psi_{u_2}(u_1,u_2,s)+(-\lambda
-s)\psi(u_1,u_2,s)}\\&& +\lambda \int_0^{u_1/\delta_1}\mu e^{-\mu
z}\psi(u_1-\delta_1z, u_2-\delta_2 z,s)dz +\lambda e^{-\mu
u_1/\delta_1}=0.
\end{eqnarray*}
In view of (\ref{eq:phiu1u2}) this equation is equivalent to
$$c_1\psi_{u_1}(u_1,u_2,s)+c_2\psi_{u_2}(u_1,u_2,s)+(-\lambda
-s)\psi(u_1,u_2,s)+\lambda \phi(u_1,u_2,s)=0.$$
This gives the first equation in the Feynman-Kac system.
To derive the second one, we apply
integration-by-parts formula to (\ref{eq:phiu1u2}):
\begin{eqnarray*}
\lefteqn{\phi(u_1,u_2,s) = \int_0^\i \frac{d}{dz}(-e^{-\mu
z})\psi(u_1-\delta_1z, u_2-\delta_2 z,s)dz}\\
&=& \psi(u_1-\d_1 z,u_2-\d_2 z,s)(-e^{-\mu
z})|_0^{\i} \\&+& \int_0^\i e^{-\mu
z}\frac{d}{dz}\psi(u_1-\delta_1z, u_2-\delta_2 z,s)dz\\
&=& \psi(u_1,u_2,s)
-\d_1\int_0^{\i}\psi_{u_1}(u_1-\delta_1z, u_2-\delta_2
z,s)e^{-\mu
z}dz\\
&-& \d_2\int_0^{\i}\psi_{u_2}(u_1-\delta_1z,
u_2-\delta_2 z,s)e^{-\mu z}dz\\
&=& \psi(u_1,u_2,s) - \d_1 \m^{-1}\phi_{u_1}(u_1, u_2,s) -
\d_2\m^{-1}\phi_{u_2}(u_1, u_2,s),
\end{eqnarray*}
which gives the second equation in the Feynman-Kac formula. The
boundary conditions follow immediately. \halmos

The above system may also be reformulated as a second
partial-differential equation in terms of $\psi(u_1,u_2,s)$ only.
To that end, we define a linear transformation $(\chi,\xi)$ of
$(\psi,\phi)$ by $\chi(r,w,s)=\psi(x,y;s)$ and $\xi(r,w,s) =
\phi(x,y;s)$ where $(x,y) = (x(r,w), y(r,w))$ are given by
$$
\left(\begin{array}{c}x \\y \end{array}\right)= \;
\left(\begin{array}{cc}\d_1&c_1\\\d_2&c_2 \end{array}\right) \;
\left(\begin{array}{c}r\\ w \end{array}\right)
$$
with inverse
transformation:
$$\left(\begin{array}{c}r\\ w \end{array}\right)
=d^{-1} \left(\begin{array}{cc}c_2&-c_1\\-\d_2&\d_1
\end{array}\right) \; \left(\begin{array}{c}x \\y
\end{array}\right)$$ where $d:=\delta_1 c_2-\delta_2 c_1 < 0$. In
the next result a PDE is derived for $\chi$.

\begin{Cor} The function
$$h(r, w):=e^{ {\m} r } e^{-(\lambda +s) w} \;  {\chi}(r,w;s)$$
solves the equation
\begin{equation}
 h_{r w}  + \m \l \;   h =0
\end{equation}
with the boundary conditions:
\begin{equation*}
\left\{\begin{array}{ll}  {h}(r,0)&= C_2e^{-(\g_2-\mu) \;
 \; r}\;
 \qquad\mbox{for all ${r }\ge 0$,}\\
h_w(r,-\frac{\d_1 r }{c_1})&= -\l e^{ \m r - (\l +s) \frac{\d_1 r
}{c_1}}\qquad\mbox{for all ${r }\geq 0$.}
\end{array}
\right. \end{equation*}
\end{Cor}
{\it Proof:} Note  that:
 \bi \im  during drift periods, $r$ is constant  and $w$
 increases, at rate $w'=1$ \im  during jump periods,  $w$ is constant  and $r$
 decreases, at rate $r'=-1$. \ei

Thus, in these new coordinates, the time is split between moving
into the direction of the axes and moving away from the axes. In
particular, at any time $T$ we have $T=T_w + T_r$, where $T_w/T_r$
are the total times of growing reserves (upward
drifting)/shrinking reserves (jumping). Note that $w=w_0 + T_w$
and $r =r_0 -T_r$ hold.

In terms of the $(r, w)$ variables, the Feynman-Kac system becomes:
\begin{equation}\label{main1}
\left(\begin{array}{l} {\chi}_{w}\\ - \xi_{r } \end{array}\right)
+\left(\begin{array}{ll}-\lambda-s & \lambda\\ \mu & -\mu
\end{array}\right) \left(\begin{array}{l} {\chi}\\ \xi
\end{array}\right)=
\left(\begin{array}{l}0\\ 0
\end{array}\right)
\end{equation}
with
\begin{equation}\label{bound1}
\left\{\begin{array}{ll}  {\chi}({r },0,s)&=
C_2e^{-\g_2 \;r },\\
\xi(r ,-\frac{\d_1}{c_1} r,s)&= 1,
\end{array}
\right.
\end{equation}
for all $r \ge 0$. Recall that the upper cone is described be
inequalities $r\ge 0$ and $w\leq 0$. Following the steps
used in the proof of differentiability of $\psi$
one can prove that $\chi$ is in class $\mathcal{C}^2$.

Eliminating $\xi$, we find \beq {\chi}_{rw} -(\l +s) {\chi}_{r} +
\mu {\chi}_{w} -s \mu \chi=0 \eeq with
\beq \left\{\begin{array}{ll}  {\chi}({r },0,s)= C_2e^{-\g_2\;r },\\
\l^{-1}\{(\l+s) \chi(r ,-\frac{\d_1}{c_1} r,s)-\chi_w(r
,-\frac{\d_1}{c_1} r)\}= 1.
\end{array}
\right. \eeq
We may remove the linear terms by switching to the
function $h$ in terms of which we get the stated result.
\halmos
\sec{Probabilistic solution \label{s:pro}}

One way to obtain the Laplace transform of the joint ruin probability
is to solve the above systems numerically. Here, we
pursue a different approach, by establishing first a general analytical
representation of the solution that holds for a general
spectrally positive L\'{e}vy process $S$.

Noting that the process $(X_1,X_2) = (U_1/\d_1,U_2/\d_2)$ has
the same ruin probability as the original two-dimensional process
$(U_1,U_2)$, we can restrict ourselves without loss of generality
to the process $(X_1,X_2)$. In the sequel we will write $\psi(x_1,x_2)$ for the
joint ruin probability (\ref{eq:ruint}) - (\ref{eq:ruin}) corresponding to
the process $(X_1,X_2)$.

\begin{Pro}\label{thseries}
If  $x_2> x_1$, then
\begin{equation}\label{eq:psibar}
\overline{\psi}(x_1,x_2):=1-\psi(x_1,x_2) =\int_0^\i
\overline{\psi}_2(z)\, \WT{P}_{(x_1,T)}(\td z),
\end{equation}
where
\begin{eqnarray}\label{T}
T &=& T(x_1, x_2) = \frac{x_2-x_1}{p_1-p_2},\\
\WT P_{(x_1, T)}(\td z) &=& P_{x_1}\left(\inf_{s\leq T} X_1(s)>0,
X_1(T)\in dz\right)\nn
\end{eqnarray}
with $P_{x}$ denoting $P$ conditioned on $\{X_1=x\}$.
\end{Pro}

{\it Proof: } In view of the definition (\ref{eq:ruint})
of $\tau$ we see that
$$
\overline{\psi}(x_1,x_2) = P\le(\min\{X_1(t), X_2(t)\} \ge
0 \mathrm{\ for\ all\ } t\ge0\ri),
$$
where $(X_1(0),X_2(0))=(x_1,x_2)$. Next, we note that,
if $x_2 > x_1$, it holds that the minimum
$$
\min\{X_1(t), X_2(t)\} =  \min\{x_1 - x_2 + (p_1-p_2)t,0\} + X_2(t)
$$
is equal to $X_1(t)$ for $t\leq  T$ and $X_2(t)$ for $t > T$,
where $T$ was defined in (\ref{T}). We have also $X_1(T)=X_2(T)$.
Subsequent application of the Markov property of $X_2$ at time $T$
shows that
\begin{eqnarray*}
\overline{\psi}(x_1,x_2) &= & P_{(x_1,x_2)}(X_1(t)\ge 0 \mathrm{\
for\ }
t\leq  T, X_2(t) \ge 0 \mathrm{\ for\ } t\ge T)\\
&= & \int_0^\i P_{x_1}\le(X_1(T) \in \td z, \inf_{s\leq T} X_1(s) \ge
0\ri) \, P_{z}\le(\inf_{s\ge 0} X_2(s) \ge 0\ri)
\\
&= & \int_0^\i P_{x_1}\le(X_1(T) \in \td z, \inf_{s\leq T} X_1(s) \ge
0\ri) \, \overline{\psi}_2(z).
\end{eqnarray*}
\halmos
\medskip

\noindent In Section \ref{s:LT} we obtain the double Laplace transform of
$\psi(x_1,x_2)$ in $x_1, x_2$, which
we invert in Section \ref{sec:LT3},
in the case of exponentially distributed jumps,
using Bromwich type contours.

\sec{\label{s:LT} Double Laplace transform in space}

Let $S(t)$ now be a general spectrally positive L\'{e}vy process and
denote by $\kappa_i(\theta)$ the Laplace exponent of the spectrally
negative L\'{e}vy process $X_i(t) = p_i t - S(t)$,
\begin{equation}\label{Lapexp}
E[e^{\theta X_i(t)}] = e^{\kappa_i(\theta)t},\qquad i=1,2.
\end{equation}
We may obtain directly the double Laplace transform in space of
$\psi(x_1,x_2)$, by exploiting
for $x_2 > x_1$ the integral representation
in Proposition \ref{thseries} and for
$x_2 \leq x_1$ the explicit formula of the Laplace
transform in $x$ of the one-dimensional ultimate
ruin probability $\psi_2(x)$. We will  use the following results:
\medskip

\begin{enumerate} \im  If $\k_i(0+)>0$, the
Laplace transform with respect to the starting point of the
ultimate survival probability $\overline{\psi}_i (x)
=P_x(\inf_{t\ge 0} X_i(t)>0)$ (see e.g. Bertoin \cite{bertoin96},
Thm. VII.8): \beq\label{LTsurv}
\la{ust} (\srui_i)^*(\th):=\int_0^\infty
e^{-\theta
x}\overline{\psi}_i(x)\;dx=\frac{\k^\prime_i(0+)}{\k_i(\th)}
\quad i=1,2.\eeq
 \im The resolvent of a spectrally negative L\'{e}vy process killed as
it enters the nonpositive half-line, due to Suprun \cite{Suprun}
(see also Bertoin \cite[Lem. 1]{bertoin97}):
\begin{eqnarray}
\lefteqn{\int_0^{\i}e^{-qt}P_{x_1}\le(\inf_{s\leq t} X_1(s)>0,
X_1(t)\in dz\ri)dt}\nonumber\\&&= \left[\exp\{-q^+(q)z\}W^{(q)}(x_1)-
{\bf 1}_{\{x_1\geq z\}}W^{(q)}(x_1-z)\right]dz, \label{resolvent}
\end{eqnarray}
where $q^+(q)$ largest root of $\kappa_1(\a)=q$ and $W^{(q)}:
[0,\i) \to [0,\i)$ is a continuous and increasing function (called
the $q$-{\it scale function} of $X_1(t)$) with the Laplace
transform:
\begin{equation}\label{eq:scale}
\I_0^\i \te{-\alpha x} W^{(q)} (y)  \td y = (\kappa_1(\alpha) -
q)^{-1},\q\q\a > q^+(q).
\end{equation}
\end{enumerate}
Now we obtain the double Laplace transform of the non-ruin probability
with respect to the initial reserves:
$$\tilde{\psi}(p,q)=\int_0^\infty\int_0^\infty e^{-px_1}
e^{-qx_2}\;\overline{\psi}(x_1,x_2)\;dx_1\;dx_2.$$
Note that
\begin{multline*}
\tilde{\psi}(p,q) =
  \I_0^\i\I_0^{x_1}e^{-px_1} e^{-qx_2}
\ovl\psi(x_1, x_2)\td x_2\td x_1 \\
+ \I_0^\i\I_{x_1}^\i e^{-px_1}
e^{-qx_2} \ovl\psi(x_1, x_2)\td x_2\td x_1.
\end{multline*}
The first Laplace transform is given by
$$
\int_0^\i\int_0^{x_1}e^{-px_1} e^{-qx_2}\srui_2 (x_2)\td x_2\td
x_1 = \frac{1}{p} (\srui_2)^*(p+q) := A.
$$
Writing $s=p+q$ and $r = (p_1-p_2)q$ we see from (\ref{eq:psibar})
and (\ref{resolvent}) that the second Laplace transform is given
by
\begin{eqnarray*}&&
\int_0^\i \int_{x_1}^\i e^{-px_1} e^{-qx_2} \ovl\psi(x_1, x_2)\td
x_2\td x_1 \nonumber\\&=& (p_1-p_2) \int_0^\i \te{-sx_1}\,dx_1
\int_0^\i \srui_2 (z) [\te{-q^+(r) z}W^{(r)}(x_1) - \mbf
1_{\{z\leq x_1\}}
W^{(r)}(x_1-z)] \td z\nonumber\\
& =& \frac{p_1-p_2}{\lape_1(s) -
r}\le[(\srui_2)^*(q^+(r))-(\srui_2)^*(s)\ri] := C - B,\label{C}
\end{eqnarray*}
where for the calculation of quantity $B$ we used (\ref{ust})
and (\ref{eq:scale}). In view of (\ref{ust}) we note that
the quantity $A-B$ is equal to \begin{equation}\label{AminB}
\frac{(\srui_2)^*(s) \k_2 (s)}{p(\lape_1(s) - r)}=
\frac{\k_2^\prime(0+)}{p(\lape_1(s) - r)}.%= \frac{\E
%X_2(1)}{p(\lape_1(s) - r)}
\end{equation}
Similarly, since $\kappa_2(\theta)=\kappa_1(\theta)+(p_2-p_1)\theta$,
we see that $C$ can be written as \bea
&&\frac{\k_2^\prime(0+)(p_1-p_2)}{\lape_2(q^+(r)) (\k_1(s) - r)}
=\frac{\k_2^\prime(0+)(p_1-p_2)}{[\lape_1(q^+(r)) +(p_2-p_1)q^+(r)
] \; (\k_1(s) - r)}\\&&=\frac{\k_2^\prime(0+)(p_1-p_2)}{(\k_1(s) -
r) (r +(p_2-p_1) q^+(r))}.
\eea
Putting everything together we find:
\begin{Pro} The double Laplace transform $\tilde{\psi}$ is given by
\begin{eqnarray}
{\tilde{\psi}(p,q)} &=&
\frac{[\k_1^\prime(0+)+(p_2-p_1)][r +(p_1-p_2)(p-q^+(r))]}{p[r
+(p_2-p_1)q^+(r) ] \; (\k_1(s) - r)} \nn \\
&=& \frac{\k_2'(0+)}{p(\k_1(p+q) - q(p_1-p_2))}
\left[1  + \frac{p}{q-q^+(q(p_1-p_2))}\right].\label{lap2}
\end{eqnarray}
\end{Pro}
\ssec{\label{sec:LT1} Exponential claims} In this section we
specialize the above result to the classical model (\ref{St})
where the jumps are exponentially distributed with parameter $\mu$
($\s_i\sim {\rm E}(\mu)$) and we write $p_i=c_i/\d_i$, $i=1,2$. In
this case the characteristic exponent of $X_i$ is given by
$$\kappa_i(\alpha)=p_i\alpha-
\frac{\lambda\alpha}{\mu+\alpha},\qquad i=1,2.$$ In particular, in
view of the form of $\k_1$ and $\k_2$ it can be verified that
$\k_2'(0+) = p_2 - \rho$ (with $\rho = \l/\m$) and $\k(p,q) = \log
E[e^{pX_1(1) + qX_2(1)}]$ is equal to
$$
\k(p,q) = \kappa_1(s) - r = \frac{p_1(z_1(q)-p)(z_2(q)-p)}{(\m+p+q)},
$$
where $r=(p_1-p_2)q$ and $s=p+q$ and
\begin{eqnarray}\label{z12} z_1(q)&=&
\frac{-(p_2 q + p_1(q+\g_1))-\sqrt{(p_2 q + p_1(q+\g_1))^2-4p_1 q
p_2(q + \g_2)}}{2p_1},
\\
z_2(q)&=&\frac{-(p_2 q + p_1(q+\g_1))+\sqrt{(p_2 q +
p_1(q+\g_1))^2-4p_1 q p_2(q + \g_2)}}{2p_1},\nonumber
\end{eqnarray}
with $\g_i = \m - \l/p_i$, $i=1,2$. For later reference we note that
\begin{eqnarray}\nn
z_1(0)&=&-\g_1\\
z_1(-\gamma_2)&=&\frac{\mu}{p_2}\left(\frac{p_2^2}{p_1}-\rho\right)^-,
\
z_2(-\gamma_2)=\frac{\mu}{p_2}\left(\frac{p_2^2}{p_1}-\rho\right)^+
\label{z1g2}
\end{eqnarray}
with $x^-=\min(x,0), x^+=\max(x,0)$.
Noting that $q^+(q(p_1-p_2))$ is the largest root of
$\kappa_1(\alpha) = q(p_1-p_2)$ and $z_2(q)$ is the largest root
of $\kappa_1(v+q) = q(p_1-p_2)$ we identify
$$
q^+(q(p_1-p_2)) = z_2(q) + q.
$$
In view of (\ref{lap2}) we thus arrive at:
\begin{Cor} If $S$ is given by (\ref{St}) with $\s_i\sim{\rm E}(\mu)$, then ${\tilde{\psi}}$ is given by
\begin{equation}
\label{secondcomp}
{\tilde{\psi}(p,q)} = \frac{(\m + p +q)(p_2 -\r)} {p
p_1( z_1(q) - p ) z_2(q)}.
\end{equation}
\end{Cor}

\sec{\label{sec:LT3} Spectral representation} In this subsection
we invert the Laplace transform (\ref{secondcomp}) of the ruin
probability $\psi$ for exponential claim sizes. To perform the
inversion we shall employ the method of residues. For an overview
of the theory of Laplace transforms and complex analysis see e.g.
Widder \cite{Widder} or Ahlfors \cite{AH}. The method of residues
leads to an explicit analytical representation of the survival
probability $\ovl\psi(x_1,x_2)$ given in the following theorem.
\begin{Thm}\label{thm:spectral}
Let $x_2>x_1$ and let $S$ is given by (\ref{St}) with
$\s_i\sim{\rm E}(\mu)$. Then it holds that
\[
\overline{\psi}(x_1,x_2)=
\begin{cases}
1  - C_1e^{-\gamma_1 x_1}
+\omega(x_1,x_2),& \text{if $\rho<\frac{p_2^2}{p_1}$,}\\
1  - C_1e^{-\gamma_1 x_1} -C_2e^{-\gamma_2 x_2}+
\frac{p_2}{p_1}e^{-\gamma_3 x_1 -\gamma_2x_2}+
\omega(x_1,x_2)&\text{else,}
\end{cases}
\]
where $\gamma_3=\frac{\mu}{p_2}\left(\rho-\frac{p_2^2}{p_1}\right)$ and
\begin{equation}\label{eq:omega}
\omega(x_1,x_2) = \frac{p_2-\rho}{\pi}\int_{q_+}^{q_-} e^{x_1a(q)+x_2 q}
\frac{f(q)\sin(b(q)x_1) + b(q)\cos(b(q)x_1)}{q(q p_2 +\mu p_2 - \l)}dq
\end{equation}
with
$$q_\pm=-\frac{1}{p_1-p_2}(\sqrt{\lambda}\pm\sqrt{p_1\mu})^2,$$
$f(q) = (\m + q + a(q))$ and
\begin{eqnarray*}
a(q)&=&\frac{-(p_1\mu-\lambda +p_2q+p_1q)}{2p_1},\\
b(q)&=&\frac{\sqrt{4p_1(p_2q\mu+p_2q^2-\lambda q)-(p_1\mu-\lambda
+p_2q+p_1q)^2}}{2p_1}.
\end{eqnarray*}
\end{Thm}
To prove this result, first observe that $\overline{\psi}(x_1,x_2)$
can be recovered from $\tilde{\psi}(p,q)$ using Mellin's formula, as folows:
\begin{equation}\label{ltdf}\overline{\psi}(x_1,x_2)=
\left(\frac{1}{2\pi
i}\right)^2\int_{\alpha-i\infty}^{\alpha+i\infty}\int_{\alpha-i\infty}^{\alpha+i\infty}
\tilde{\psi}(p,q)e^{x_2 q} e^{x_1 p}\;dp\;dq,
\end{equation}
where $\alpha>0$.
The next step consists in iteratively evaluating this double
integral (first w.r.t. $p$ and then w.r.t. $q$) using Cauchy's theorem.
The result of the first inversion is given in the next result:
\begin{Lem}\label{lem:psis4}
For $\a>0$ and fixed $q$ with $\Re(q)>0$ it holds that
\begin{equation}\label{eq:inv1}
\frac{1}{2\pi i}\I_{\a - i\i}^{\a + i\i}\T\psi(p,q)e^{x_1p}dp =
\frac{\k_2'(0+)}{\k_2(q)} - e^{x_1z_1(q)}g(q),
\end{equation}
where $g(q)$ is given by \begin{equation}\label{gq}
g(q)= \frac{(p_2-\rho)(\m + z_1(q) + q)}{q(\m p_2 - \l + p_2q)}.
\end{equation}
\end{Lem}
The proofs of results that are not developed in the text can be
found in the Appendix.

In view of equations (\ref{LTsurv}) and (\ref{exponruin}) we
recognize the first term in (\ref{eq:inv1}) as the Laplace
transform of $\ovl\psi_2(x) = 1 - C_2\te{-\g_2 x}$. The inversion
of the second term relies on the following properties of $g(q)$
and $z_1(q)$, that were defined in (\ref{gq}) and (\ref{z12})
respectively:

\begin{Lem}\label{lem:psis3}
(i) The functions $g$ and $z_1$,
are analytic in the set\
$$
Q=\{q\in\mathbb C: q\notin\{ 0,-\g_2\}\cup[q_+,q_-] \},
$$
where  $-\g_2 > q_-$ and $qg(q)$ remains bounded if $|q|\to\i$.

(ii) Let $q^\pm_\e = q \pm i\e$ with $q\in[q_+,q_-]$ and $\e>0$.
If $\e\downarrow0$, then
\begin{equation}\label{z1qe}
z_1(q^+_\e) \to z^-(q):=a(q)-ib(q),\q z_1(q^-_\e) \to
z^+(q):=a(q)+ib(q).
\end{equation}
\end{Lem}
\begin{figure}[t]
\begin{center}
\leavevmode \resizebox{0.4
\textwidth}{!}{\includegraphics{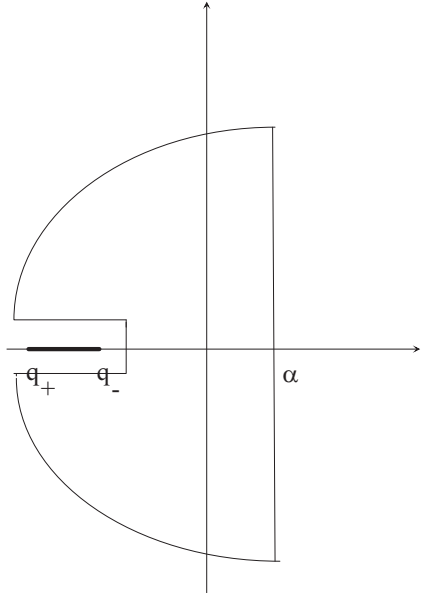}} \caption{Bromwich
contour}\label{fig2}
\end{center}
\end{figure}
In order to ensure that we can calculate the inversion
of the second term in (\ref{eq:inv1}) using the method of residues
we fix $0<a<\g_2$ and replace $g(q)$ by $g(q)/(q+a)$ (note that
in view of Lemma \ref{lem:psis3} the latter is
$O(q^{-2})$ as  $|q|\to\i$).
Denote by $f$ and $h_a$ the Laplace inverses of $g(q)$
and $g(q)/(q+a)$, that is
$$
g(q)/(q+a) = \int_0^\i e^{-qx}h_a(x)dx =
\int_0^\i e^{-qx}\int_0^x e^{-a(x-y)}f(y)dydx.
$$
Then it follows that $f$ can be recovered from $h_a$ by
\begin{equation}\label{eq:fh}
f(x) = \lim_{a\downarrow 0}\frac{d}{dx}h_a(x).
\end{equation}
To complete the inversion of $\T\psi(p,q)$ we are thus led to
evaluate the integral
\begin{equation}\label{eq:Brint}
\frac{1}{2\pi i}\I_{\a - i\i}^{\a + i\i}k_a(q)dq
\q\text{where}\q k_a(q):=\frac{g(q)}{q+a}e^{x_2q+x_1z_1(q)}.
\end{equation}
Using Lemma \ref{lem:psis3} we choose a Bromwich contour $
\Gamma_{R,\e}$ that encloses the poles of $e^{z_1(q)x_1}g(q)$
while the cut of the square root $z_1(q)$ is not enclosed (see
Figure \ref{fig2} and the proof of Lemma \ref{lem:psis6} for
formal definition of $\Gamma_{R,\e}$). We note that this is a
standard approach to calculate integrals of the form
(\ref{eq:Brint}) (e.g. Ahlfors \cite{AH} or see Pervozvansky
\cite{Per} for a recent application to the calculation of
one-dimensional ruin probabilities).

Recalling from Lemma \ref{lem:psis3} that $g(q)$ has two simple
poles, in $q=0$ and $q=-\g_2$, Cauchy's theorem implies
that
\begin{equation}\label{eq:CAU}
\frac{1}{2\pi i}\oint_{\Gamma_{R,\e}}k_a(q)dq =
\text{Res}_{q=0}k_a(q) + \text{Res}_{q=-\g_2}k_a(q) +
\text{Res}_{q=-a}k_a(q).
\end{equation}
The next step consists in evaluating the residues in (\ref{eq:CAU}),
which is a matter of straightforward calulations:

\begin{Lem}\label{lem:psis5} Writing
$\T C_2 = C_2 + \frac{z_1(-\g_2)}{\m}$ it holds for $a>0$ that
\begin{eqnarray*}
\label{eq:residue1}
\mathrm{Res}_{q=-a}k_a(q) &=& g(-a)e^{-x_2a+x_1z_1(-a)},\
\mathrm{Res}_{q=0}k_a(q) =  -\frac{\l}{a\m p_1}\te{-\g_1x_1}\\
\label{eq:residue3} \mathrm{Res}_{q=-\g_2}k_a(q) &=&
\frac{\T C_2}{a-\g_2}\te{z_1(-\g_2)x_1-\g_2x_2}.
\end{eqnarray*}
\end{Lem}
Next we turn to the left-hand side of the formula (\ref{eq:CAU}):

\begin{Lem}\label{lem:psis6}
For $\a>0$ it holds that
\begin{equation*}
\lim_{\e\downarrow 0}\lim_{R\to\i}\frac{1}{2\pi
i}\oint_{\Gamma_{R,\e}}k_a(q)dq = \frac{1}{2\pi i}\I_{\a -
i\i}^{\a + i\i}k_a(q)dq - \omega_a(x_1,x_2),
\end{equation*}
where  $\omega_a(x_1,x_2)$ is given by (\ref{eq:omega}) but with
$dq$  replaced by $dq/(q+a)$.
\end{Lem}

{\it Proof of Theorem \ref{thm:spectral}}
In view of (\ref{eq:fh}), (\ref{eq:CAU}) and Lemmas
\ref{lem:psis5} and \ref{lem:psis6}, the final result is obtained
by first differentiating the residues in Lemma \ref{lem:psis5} and $\omega_a$
with respect to $x_2$ and subsequently letting subsequently $a$ tend to zero.
Taking note of the facts
$$
\lim_{a\downarrow 0}z_1(-a) = z_1(0) = -\g_1,\
\lim_{a\downarrow 0}ag(a)=\WT C_2\ \text{and}\
\lim_{a\downarrow 0}\frac{\partial}{\partial x_2}
\w_a(x_1,x_2) = \omega(x_1,x_2)
$$
 completes the proof 
(where the latter follows using the dominated convergence thoerem). \halmos

\appendix
\section{\label{sec:app} Appendix}
\subsection{Formal construction of the fluid-embedding}
A formal construction of the process $(\T U_1, \T U_2)$ is as follows.
Let $\tau_1, \tau_2, \ldots $ and $\s_1,\s_2,\ldots$ denote the subsequent
inter-arrival times and claim sizes. Note that these form
sequences of i.i.d. exponential random variables with parameters
$\l$ and $\mu$ respectively. Define the switch times $S_n$ by 
$S_0=0$ and, for $n\ge 1$,
\begin{eqnarray*}
S_{2n-1} &=& S_{2n-2} + \tau_n, S_{2n} = S_{2n-1} + \sigma_n \quad\text{if $Y_0=1$}\\
S_{2n-1} &=& S_{2n-2} + \sigma_n, S_{2n} = S_{2n-1} + \tau_n \quad\text{if $Y_0=-1$}
\end{eqnarray*}
%For $n\ge 1$ define the $n$th switch
%time $S_n$ by
%$$
%S_{2n-1} = S_{2n-2} + \t_{n}\q \text{ and }\q
%S_{2n} = S_{2n-1} + \s_{n}\q\text{with $S_0=0$}
%$$
and construct $Y$ taking values in $\{-1,+1\}$ by setting
$$
Y_{S_n} = - Y_{S_{n-1}}
\q\text{for $n\ge 1$ with $Y_0\in\{-1,+1\}$}.
$$
Then $Y$ is a two-state
Markov chain indicating whether $(\T U_1, \T U_2)$ is
increasing (state $+1$) or decreasing (state $-1$); more precisely,
denoting by
$$I(t) = \I_0^t\mbf 1_{\{Y_s=1\}}\td s $$ the total time up to $t$ that
$Y$ has spent in state $+1$, we set
$$\T U_i(t) = p_i I(t) - \d_i(t - I(t)),\q i=1,2,$$
and the construction is complete.

\subsection{Proof of Lemma \ref{lem:psis4}}
By performing a partial fraction decomposition (in $p$) it follows that
\begin{eqnarray}\nn
\T\psi(p,q) &=& \frac{1}{p}\frac{(p_2-\rho)(\m+q)}{q(\m p_2 - \l + qp_2)}
- \frac{1}{p-z_1(q)}\frac{(p_2-\rho)(\m+q+z_1(q))}{q(\m p_2 - \l + qp_2)}\\
&=& \frac{1}{p} \frac{\k_2'(0)}{\k_2(q)} - \frac{1}{p-z_1(q)}g(q).
\label{1op}
\end{eqnarray}
Since $\I_0^\i e^{-pt}e^{ct}dt=(p-c)^{-1}$ for $p>c$,
the result follows by inverting (\ref{1op}) term by term. \halmos

\subsection{Proof of Lemma \ref{lem:psis3}}
(i) Noting that the argument of the square root in (\ref{z12}) is positive
if $q$ is real and $q < q_+$ or $q> q_-$, it follows that $z_1(q)$
is analytic outside the cut $[q_+,q_-]$. Since, furthermore,
the denominator of $g$ has no roots in $Q$, we see that $g(q)$ is analytic
in the set $Q$. The asymptotics directly follow from the form of $g$.

(ii) Employing the standard
definition of the square root $z\mapsto\sqrt{z}$
(with the cut along the negative half-line) and
appealing to the definition of $z_1$ and the continuity of
the argument $\text{Arg}(z)$ and modulus $|z|$
imply that the convergence in (\ref{z1qe}) holds true. \halmos

\subsection{Proof of Lemma \ref{lem:psis6}}
Consider the contour $\G_{R,\e}$ that is given in Figure
\ref{fig2}, i.e. $\G_{R,\e}$ consists of the line segments $[\a-i
R,\a + iR]$, $[q_-+\e -i \e, q_-+\e +i\e]$, $[q_+-\e-i\e,q_-+\e
-i\e]$ and $[q_+-\e+i\e,q_-+\e +i\e]$ and two quarter circles in
the left half-plane joining $q_+-\e + i\e$ and $\a + iR$, and,
$q_+-\e - i\e$ and $\a - iR$, respectively.

By taking the limits of $R\to\i$ and subsequently letting
$\e\downarrow 0$ and using that the integrals of the
quarter-circles tend to zero (in view of the fact that
$g(q)/(q+a)=O(q^{-2})$ as $|q|\to\i$, cf. Lemma
\ref{lem:psis3}(i))  we find that the contour integral
$\frac{1}{2\pi i}\oint_{\Gamma_{R,\e}} k_a(q)\;dq$ converges to
\begin{equation}\label{eq:3}
\frac{1}{2\pi i}\int_{\a-i\i}^{\a + i\i}k_a(q)\;dq + \frac{1}{2\pi
i}\int_{q_+\to q_-} k_a(q)\;dq + \frac{1}{2\pi i}\int_{q_-\to q_+}
k_a(q)\;dq,
\end{equation}
where integrals $\int_{q_+\to q_-}$ and $\int_{q_-\to q_+}$ are
the limits of the line integrals along the segments
$[q_+-\e+i\e,q_-+\e +i\e]$ and $[q_--\e-i\e,q_++\e -i\e]$ of the
contour $\G_{R,\e}$, respectively. In view of Lemma
\ref{lem:psis3}(ii) it follows that the last two integrals in
(\ref{eq:3}) are equal to
\begin{multline}
 \frac{p_2-\rho}{2\pi i}\int_{q_+}^{q_-}\frac{\m+q+a(q)}{q(q+a)(p_2q + \m p_2 - \l)}
e^{x_2q}(e^{z^-(q)x_1} - e^{z^+(q)x_1})dq\\
+ \frac{p_2-\rho}{2\pi }\int_{q_+}^{q_-}\frac{b(q)}{q(q+a)(p_2q +
\m p_2 - \l)} e^{x_2q}(e^{z^-(q)x_1} +
e^{z^+(q)x_1})dq,\label{iiint}
\end{multline}
where $z^-(q)$ and $z^+(q)$ are defined in (\ref{z1qe}). Use of
the representation $e^{a+ib} = e^a(\cos b + i\sin b)$ (for
$a,b\in\mathbb R$) completes the calculation of the contour
integral of $g(q)/(q+a)$.
 \halmos

\subsection*{\bf Acknowledgements}

\no F. Avram and M. Pistorius gratefully acknowledge support from
the London Mathematical Society, grant \# 4416. F. Avram and Z.
Palmowski acknowledge support by POLONIUM no 09158SD. Z. Palmowski
acknowledges support by KBN 1P03A03128 and NWO 613.000.310. M.
Pistorius acknowledges support by the Nuffield Foundation
NUF/NAL/000761/G.

\end{document}